\documentclass{article}
\usepackage{amsmath}
\usepackage{amsfonts}
\usepackage[ansinew]{inputenc}
\usepackage{graphicx}

\usepackage{harvard}

\numberwithin{equation}{section}

\begin{document}

\title{Recursive Numerical Evaluation of the Cumulative Bivariate Normal Distribution}
\author{Christian Meyer\thanks{The author wishes to thank Axel Vogt for helpful discussion.} \footnote{DZ BANK AG, Platz der Republik, D-60265 Frankfurt. The opinions or recommendations expressed in this article are those of the author and are not representative of DZ BANK AG.} \footnote{E-Mail: {\tt Christian.Meyer@dzbank.de}}}
\date{\today}
\maketitle

\begin{abstract}
We propose an algorithm for evaluation of the cumulative bivariate normal distribution, building upon Marsaglia's ideas for evaluation of the cumulative univariate normal distribution. The algorithm is mathematically transparent, delivers competitive performance and can easily be extended to arbitrary precision.
\end{abstract}

\section{Introduction}

The cumulative normal distribution, be it univariate or multivariate, has to be evaluated numerically. There are numerous algorithms available, many of these having been fine-tuned, leading to faster evaluation and higher accuracy but also to lack of mathematical transparency.

For the univariate case, \citeasnoun{Marsaglia} has proposed a very simple and intuitive but powerful alternative that is based on Taylor expansion of Mills' ratio or similar functions. In this note we will extend Marsaglia's approach to the bivariate case. This will require two steps: reduction of the evaluation of the cumulative bivariate normal distribution to evaluation(s) of a univariate function, i.e., to the cumulative bivariate normal distribution on the diagonal, and Taylor expansion of that function. Note that a similar approach, but with reduction to the axes instead of the diagonals, has been proposed by \citeasnoun{Vogt}.

The resulting algorithm has to be compared with existing approaches. For overview on and discussion of the latter, cf. \cite{BretzGenz}, \cite{AC}, \cite{TW}, and \cite{WangKen}. Most implementations today will rely on variants of the approaches of \citeasnoun{Divgi} or of \citeasnoun{DW}. Improvements of the latter method have been provided by \citeasnoun{Genz} and \citeasnoun{West}. The method of \citeasnoun{Drezner}, although less reliable, is also very common, mainly because it is featured in \cite{Hull} and other prevalent books.

It will turn out that the algorithm proposed in this paper is able to deliver near double precision (in terms of absolute error) using double arithmetic. Furthermore, implementation of the algorithm using high-precision libraries is straightforward; indeed, a quad-double implementation has been applied for testing purposes. Performance is competitive, and trade-offs between speed and accuracy may be implemented with little effort.


\section{Theory}
\label{sec_theory}

In this section we are going to develop the algorithm. In order to keep the presentation lean we will often refer to the author's recent survey \cite{Meyer}. For further background on normal distributions the reader is also referred to text books such as \cite{BL}, \cite{KBJ} and \cite{PR}.

\subsection{Evaluation on the diagonal}
\label{subsec_diagonal}

Denote by
\[
\varphi(x) := \frac{1}{\sqrt{2\pi}} \exp\left(-\frac{x^2}{2}\right),
\qquad
\Phi(x) := \int_{-\infty}^x \varphi(t)\;dt
\]
the density and distribution function of the standard normal distribution. Mills' ratio is then defined as
\[
R(x) := \frac{1-\Phi(x)}{\varphi(x)} = \frac{\Phi(-x)}{\varphi(-x)}.
\]
Furthermore, denote by
\begin{align*}
\varphi_2(x,y;\varrho) & := \frac{1}{2\pi\sqrt{1-\varrho^2}}
\exp\left(-\frac{x^2-2\varrho xy + y^2}{2(1-\varrho^2)}\right),\\
\Phi_2(x,y;\varrho) & := \int_{-\infty}^{x} \int_{-\infty}^{y}
\varphi_2(s,t;\varrho) \; dt \; ds,
\end{align*}
the density and distribution function of the bivariate standard normal distribution with correlation parameter $\varrho\in(-1,1)$. We will also write
\begin{align*}
\varphi_2(x;\varrho) & := \varphi_2(x,x;\varrho) = \frac{1}{2\pi\sqrt{1-\varrho^2}}\exp\left(-\frac{x^2}{1+\varrho}\right),\\
\Phi_2(x;\varrho) & := \Phi_2(x,x;\varrho),\\
\lambda(\varrho) & := \sqrt{\frac{1-\varrho}{1+\varrho}}.
\end{align*}
We are going to use the following properties:
\begin{align*}
\frac{d}{dx} \varphi_2(x;\varrho) & = -\frac{2x}{1+\varrho}\cdot \varphi_2(x;\varrho),\\
\frac{d}{dx} \Phi_2(x;\varrho) & = 2\cdot\varphi(x)\cdot\Phi(\lambda(\varrho)\cdot x),\\
\varphi_2(x;\varrho) \cdot\sqrt{1-\varrho^2}& = \varphi(x)\cdot\varphi(\lambda(\varrho)\cdot x)
\end{align*}
In the following we will assume that $x\leq 0$, $\varrho\geq 0$. In this case the following bounds apply (cf. \cite[Th. 5.2]{Meyer}):
\begin{equation}
1+\frac{2}{\pi}\arcsin(\varrho)
\leq \frac{\Phi_2(x;\varrho)}{\Phi(x)\cdot\Phi(\lambda(\varrho)\cdot x)}\leq
1+\varrho
\label{eq_bounds}
\end{equation}
Furthermore, as is proven implicitly in \cite[App. A.2]{Meyer},
\[
\lim_{x\longrightarrow-\infty} \frac{\Phi_2(x;\varrho)}{\Phi(x)\cdot\Phi(\lambda(\varrho)\cdot x)} = 1+\varrho.
\]
Now we define
\[
D(x) := \frac{(1+\varrho)\cdot\Phi(x)\cdot\Phi(\lambda(\varrho)\cdot x)-\Phi_2(x;\varrho)}{\varphi_2(x;\varrho)}.
\]
Starting with
\[
D'(x) = (\varrho-1)\cdot\sqrt{1-\varrho^2}\cdot R(-\lambda(\varrho)\cdot x)
+(1-\varrho^2)\cdot R(-x) + \frac{2x}{1+\varrho}\cdot D(x).
\]
we find the recursion
\begin{align*}
D^{(k)}(x) & = (\varrho-1)\cdot\sqrt{1-\varrho^2}\cdot (-\lambda(\varrho))^{k-1}\cdot R^{(k-1)}(-\lambda(\varrho)\cdot x)\\
& + (1-\varrho^2)\cdot (-1)^{k-1}\cdot R^{(k-1)}(-x)\\
& + \frac{2(k-1)}{1+\varrho}\cdot D^{(k-2)}(x) + \frac{2x}{1+\varrho}\cdot D^{(k-1)}(x)
\end{align*}
which we can use to recursively evaluate the Taylor expansion of $D$ around zero. Dividing by $\sqrt{1-\varrho^2}$ for convenience, we define
\begin{align*}
a_k & := \frac{x^{k+1}}{k!}\cdot(\varrho-1)\cdot(-\lambda(\varrho))^k \cdot R^{(k)}(0),\\
b_k & := \frac{x^{k+1}}{k!}\cdot\sqrt{1-\varrho^2}\cdot(-1)^k \cdot R^{(k)}(0),\\
d_k & := \frac{x^k}{k!}\cdot\frac{D^{(k)}(0)}{\sqrt{1-\varrho^2}}.
\end{align*}
Using
\[
R^{(k)}(x) = (k-1)\cdot R^{(k-2)}(x) + x\cdot R^{(k-1)}(x)
\]
we derive the following recursion scheme:
\begin{align}
a_k & = \frac{1}{k}\cdot x^2\cdot\frac{1-\varrho}{1+\varrho}\cdot a_{k-2},\label{eq_rec_1}\\
b_k & = \frac{1}{k}\cdot x^2\cdot b_{k-2},\label{eq_rec_2}\\
d_k & = \frac{1}{k}\cdot\left(a_{k-1}+b_{k-1}+\frac{2x^2}{1+\varrho}\cdot d_{k-2}\right),\label{eq_rec_3}
\end{align}
with initial values
\begin{align}
a_0 & = (\varrho-1)\cdot \sqrt{\frac{\pi}{2}}\cdot x,\label{eq_start_1}\\
a_1 & = \lambda(\varrho)\cdot(\varrho-1)\cdot x^2,\label{eq_start_2}\\
b_0 & = \sqrt{1-\varrho^2}\cdot \sqrt{\frac{\pi}{2}}\cdot x,\label{eq_start_3}\\
b_1 & = \sqrt{1-\varrho^2}\cdot x^2,\label{eq_start_4}\\
d_0 & = \frac{\varrho\cdot \pi}{2}-\arcsin(\varrho),\label{eq_start_5}\\
d_1 & = \left(\varrho-1+\sqrt{1-\varrho^2}\right)\cdot\sqrt{\frac{\pi}{2}}\cdot x.\label{eq_start_6}
\end{align}
Here we have used that
\[
R(0)=\sqrt{\frac{\pi}{2}}, \qquad \Phi_2(0;\varrho)=\frac{1}{4}+\frac{1}{2\pi}\cdot\arcsin(\varrho).
\]
We can now compute $\Phi_2(x;\varrho)$ numerically via
\[
\Phi_2(x;\varrho) = (1+\varrho)\cdot \Phi(x)\cdot\Phi(\lambda(\varrho)\cdot x) - \frac{1}{2\pi}\cdot \exp\left(-\frac{x^2}{1+\varrho}\right)\cdot\left(\sum_{k=0}^{\infty} d_k\right).
\]
Note that it would also have been possible to work with, e.g., one of the functions
\begin{align*}
D_2(x) & := \frac{1-\Phi_2(x;\varrho)}{\varphi_2(x;\varrho)},\\
D_3(x) & := \frac{\Phi_2(0;\varrho)-\Phi_2(x;\varrho)}{\varphi_2(x;\varrho)},\\
D_4(x) & := \frac{2\cdot\Phi(x)\cdot\Phi(\lambda(\varrho)\cdot x)-\Phi_2(x;\varrho)}{\varphi_2(x;\varrho)}
\end{align*}
instead. The resulting recursion schemes are in fact easier (two summands instead of three) but will be running into numerical problems (cancellation, or lower accuracy for $x\longrightarrow -\infty$).

\subsection{Reduction to the diagonal}
\label{subsec_reduction}

In order to apply the results from Section \ref{subsec_diagonal} to the numerical evaluation of $\Phi_2(x,y;\varrho)$ for general $x$, $y$ and $\varrho$, we start with the symmetric formula (cf. \cite[Eq. (3.16)]{Meyer})
\begin{equation}
\Phi_2(x,y;\varrho) = \Phi_2(x,0;\varrho_x) - \delta_x + \Phi_2(0,y;\varrho_y) - \delta_y,
\label{eq_2axis}
\end{equation}
where
\[
\delta_x =
\begin{cases}
\frac{1}{2}, & x < 0 \quad \text{and}\quad y\geq 0,\\
0, & \text{else},
\end{cases}, \qquad
\delta_y =
\begin{cases}
\frac{1}{2}, & y < 0 \quad \text{and}\quad x\geq 0,\\
0, & \text{else},
\end{cases}
\]
and
\begin{align*}
\varrho_x & = -\frac{\alpha_x}{\sqrt{1+\alpha_x^2}},
\qquad \alpha_x = \frac{1}{\sqrt{1-\varrho^2}}\left(\frac{y}{x}-\varrho\right),\\
\varrho_y & = -\frac{\alpha_y}{\sqrt{1+\alpha_y^2}},
\qquad \alpha_y = \frac{1}{\sqrt{1-\varrho^2}}\left(\frac{x}{y}-\varrho\right).
\end{align*}
From the axis $y=0$ to the diagonal $x=y$ we get by applying the formula (cf. \cite[Eq. (3.18)]{Meyer})
\begin{equation}
\Phi_2(x,0;\varrho) =
\begin{cases}
\frac{1}{2} \cdot\Phi_2(x,x;1-2\varrho^2), & \quad \varrho < 0,\\
\Phi(x) - \frac{1}{2} \cdot\Phi_2(x,x;1-2\varrho^2), & \quad \varrho \geq 0.
\end{cases}
\end{equation}
Specifically, we obtain
\begin{equation}
1-2\varrho_x^2 = 1-\frac{2\cdot(\varrho x - y)^2}{x^2+y^2-2\varrho xy}
=1-\frac{2\cdot a_x}{1+a_x}
\label{eq_rho_1}
\end{equation}
with
\begin{align}
a_x & = \left(\frac{\varrho x - y}{x\cdot\sqrt{1-\varrho^2}}\right)^2\\
& = \left(\frac{x-y}{x\cdot\sqrt{1-\varrho^2}}-\sqrt{\frac{1-\varrho}{1+\varrho}}\right)^2\label{eq_ax_rho1}\\
& = \left(\frac{x+y}{x\cdot\sqrt{1-\varrho^2}}-\sqrt{\frac{1+\varrho}{1-\varrho}}\right)^2\label{eq_ax_rhom1}
\end{align}
where in an implementation (\ref{eq_ax_rho1}) should be used for $\varrho\longrightarrow 1$, and (\ref{eq_ax_rhom1}) for $\varrho\longrightarrow -1$, in order to avoid catastrophic cancellation. Note also that
\[
\varrho_x < 0 \quad\Longleftrightarrow\quad \alpha_x > 0 \quad\Longleftrightarrow\quad \frac{y}{x} > \varrho.
\]
In a last step, if necessary to ensure $x\leq 0$ and $\varrho\geq 0$, we apply the formulas (cf. \cite[Eq. (2.15)]{Meyer} and \cite[Eq. (3.27)]{Meyer})
\begin{align}
\Phi_2(x,x;\varrho) & = 2\cdot \Phi(x) - 1 + \Phi_2(-x,-x;\varrho),\label{eq_minusx}\\
\Phi_2(x,x;\varrho) & = 2\cdot\Phi(x)\cdot\Phi\left(\lambda(\varrho)\cdot x\right) - 
\Phi_2\left(\lambda(\varrho)\cdot x,\lambda(\varrho)\cdot x,-\varrho\right).\label{eq_minusrho}
\end{align}
Specifically, we obtain
\begin{align}
|\lambda(1-2\varrho_x^2)\cdot x| = |x|\cdot\sqrt{\frac{\varrho_x^2}{1-\varrho_x^2}} & = \frac{|\varrho x - y|}{\sqrt{1-\varrho^2}}\\
& = \left|\frac{x-y}{\sqrt{1-\varrho^2}}-x\cdot\sqrt{\frac{1-\varrho}{1+\varrho}}\right|\label{eq_bx_rho1}\\
& = \left|\frac{x+y}{\sqrt{1-\varrho^2}}-x\cdot\sqrt{\frac{1+\varrho}{1-\varrho}}\right|\label{eq_bx_rhom1},
\end{align}
where in an implementation (\ref{eq_bx_rho1}) should be used for $\varrho\longrightarrow 1$, and (\ref{eq_bx_rhom1}) for $\varrho\longrightarrow -1$, in order to avoid catastrophic cancellation.

It will be favorable to work with
\[
2\varrho_x^2 = \frac{2\cdot a_x}{1+a_x}
\]
instead of $1-2\varrho_x^2$. If (\ref{eq_minusrho}) has to be applied (i.e., if $1-2\varrho_x^2<0$, which is equivalent with $a_x>1$), correspondingly we will work with
\begin{equation}
1-(-(1-2\varrho_x^2)) = 2-2\varrho_x^2 = \frac{2}{1+a_x}.
\label{eq_rho_2}
\end{equation}

\section{Implementation}
\label{sec_implement}

In the following we will discuss implementation of the algorithm derived in Section \ref{sec_theory}. The C++ language has been chosen because it is the market standard in quantitative finance, one of the fields frequently requiring evaluation of normal distributions.

\subsection{Evaluation on the diagonal}

Source code (in C++) for evaluation of $\Phi_2(x;\varrho)$ as in Section \ref{subsec_diagonal}, for $x\leq 0$ and $\varrho\geq 0$, is provided in Figure \ref{source_diag}. In the following we will comment on some details of the implementation.

Equations (\ref{eq_rec_1}) - (\ref{eq_start_6}) show that it is reasonable to provide $a:=1-\varrho$, instead of $\varrho$, as input for the evaluation of $\Phi_2(x;\varrho)$. Moreover, cf. (\ref{eq_rho_1}) and (\ref{eq_rho_2}), double inversion (i.e., computation of $1-(1-z)$ instead of $z$) is to be avoided in the reduction algorithm.

Values for $\Phi(x)$ and for $\Phi(\lambda(\varrho)\cdot x)$ are also expected as input parameters. This makes sense because the values are needed by the reduction algorithm as well (and hence should not be computed twice).

Evaluation of $\arcsin(\varrho)$ is to be avoided for $\varrho \longrightarrow 1$ and has been replaced (without optimization of the cutoff point) by
\[
\arcsin(\varrho) = \arccos\left(\sqrt{\frac{1-\varrho}{1+\varrho}}\right) = \arccos(\lambda(\varrho)).
\]
Note that $\lambda(\varrho)$ has to be computed anyway.

Constants (all involving $\pi$) have been pre-computed in double precision. The recursion stops if a new term does not change the computed sum. If the a priori bound for the absolute error, given by (\ref{eq_bounds}), is less than $5\cdot 10^{-17}$, the upper bound is returned (relative accuracy on the diagonal may be increased by dropping this condition but overall relative accuracy will still be determined by the reduction to the diagonal, cf. Section \ref{subsec_reduction_imp}), and by the accuracy of the implementation of $\Phi$. The final result is always checked against the upper and lower bound.

Note that $d_{2k}$ and $d_{2k+1}$ have different sign but comparable order. Bracketing them before summation can therefore reduce cancellation error.

\begin{figure}
\hrule
\vspace*{2mm}
\begin{small}
\begin{verbatim}
double Phi2diag( const double& x,
                 const double& a,      // 1 - rho
                 const double& px,     // Phi( x )
                 const double& pxs )   // Phi( lambda( rho ) * x )
{
    if( a <= 0.0 ) return px;        // rho == 1
    if( a >= 1.0 ) return px * px;   // rho == 0

    double b = 2.0 - a, sqrt_ab = sqrt( a * b );
    double asr = ( a > 0.1 ? asin( 1.0 - a ) : acos( sqrt_ab ) );
    double comp = px * pxs;
    if( comp * ( 1.0 - a - 6.36619772367581343e-001 * asr ) < 5e-17 )
        return b * comp;

    double tmp = 1.25331413731550025 * x;
    double a_coeff = a * x * x / b;
    double a_even = -tmp * a;
    double a_odd = -sqrt_ab * a_coeff;
    double b_coeff = x * x;
    double b_even = tmp * sqrt_ab;
    double b_odd = sqrt_ab * b_coeff;
    double d_coeff = 2.0 * x * x / b;
    double d_even = ( 1.0 - a ) * 1.57079632679489662 - asr;
    double d_odd = tmp * ( sqrt_ab - a );

    double res = 0.0, res_new = d_even + d_odd;
    int k = 2;
    while( res != res_new )
    {
        d_even = ( a_odd + b_odd + d_coeff * d_even ) / k;
        a_even *= a_coeff / k;
        b_even *= b_coeff / k;
        k++;
        a_odd *= a_coeff / k;
        b_odd *= b_coeff / k;
        d_odd = ( a_even + b_even + d_coeff * d_odd ) / k;
        k++;
        res = res_new;
        res_new += d_even + d_odd;
    }
    res *= exp( -x * x / b ) * 1.591549430918953358e-001;
    return max( ( 1.0 + 6.36619772367581343e-001 * asr ) * comp,
                b * comp - max( 0.0, res ) );
}
\end{verbatim}
\end{small}
\hrule
\caption{C++ source code for evaluation of $\Phi_2(x;\varrho)$}
\label{source_diag}
\end{figure}

\subsection{Reduction to the diagonal}
\label{subsec_reduction_imp}

Source code (in C++) for evaluation of $\Phi_2(x;\varrho)$ as in Equation (\ref{eq_2axis}) is provided in Figure \ref{source_Phi}, and source code for evaluation of $\Phi_2(x,0;\varrho_x) - \delta_x$ is provided in Figure \ref{source_help}. In the following we will comment on some details of the implementation.

The special cases $|\varrho|=1$ and $x=y=0$ are dealt with in {\tt Phi2()}. Therefore, in {\tt Phi2help()} there is no check against {\tt 1.0 - rho == 0.0}, {\tt 1.0 + rho == 0.0} or {\tt s == 0.0}.

It is assumed that {\tt sqr(x)} evaluates {\tt x*x}. The cutoff points $|\varrho|=0.99$ have been set by visual inspection and might be optimized.

\begin{figure}
\hrule
\vspace*{2mm}
\begin{small}
\begin{verbatim}
double Phi2help( const double& x,
                 const double& y,
                 const double& rho )
{
    if( x == 0.0 ) return ( y >= 0.0 ? 0.0 : 0.5 );

    double s = sqrt( ( 1.0 - rho ) * ( 1.0 + rho ) );

    double a = 0.0, b1 = -fabs( x ), b2 = 0.0;
    if( rho > 0.99 )
    {
        double tmp = sqrt( ( 1.0 - rho ) / ( 1.0 + rho ) );
        b2 = -fabs( ( x - y ) / s - x * tmp );
        a = sqr( ( x - y ) / x / s - tmp );
    }
    else if( rho < -0.99 )
    {
        double tmp = sqrt( ( 1.0 + rho ) / ( 1.0 - rho ) );
        b2 = -fabs( ( x + y ) / s - x * tmp );
        a = sqr( ( x + y ) / x / s - tmp );
    }
    else
    {
        b2 = -fabs( rho * x - y ) / s;
        a = sqr( b2 / x );
    }

    double p1 = Phi( b1 ), p2 = Phi( b2 );   // cum. standard normal

    double q = 0.0;
    if( a <= 1.0 )
        q = 0.5 * Phi2diag( b1, 2.0 * a / ( 1.0 + a ), p1, p2 );
    else
        q = p1 * p2 - 0.5 * Phi2diag( b2, 2.0 / ( 1.0 + a ), p2, p1 );

    int c1 = ( y / x >= rho );
    int c2 = ( x < 0.0 );
    int c3 = c2 && ( y >= 0.0 );
    return ( c1 && c3 ? q - 0.5
                      : c1 && c2 ? q
                      : c1 ? 0.5 - p1 + q
                      : c3 ? p1 - q - 0.5
                      : c2 ? p1 - q
                      : 0.5 - q );
}
\end{verbatim}
\end{small}
\hrule
\caption{C++ source code for evaluation of $\Phi_2(x,0;\varrho_x) - \delta_x$}
\label{source_help}
\end{figure}

\begin{figure}
\hrule
\vspace*{2mm}
\begin{small}
\begin{verbatim}
double Phi2( const double& x,
             const double& y,
             const double& rho )
{
    if( ( 1.0 - rho ) * ( 1.0 + rho ) <= 0.0 )   // |rho| == 1
        if( rho > 0.0 )
            return Phi( min( x, y ) );
        else
            return max( 0.0, min( 1.0, Phi( x ) + Phi( y ) - 1.0 ) );

    if( x == 0.0 && y == 0.0 )
        if( rho > 0.0 )
            return Phi2diag( 0.0, 1.0 - rho, 0.5, 0.5 );
        else
            return 0.5 - Phi2diag( 0.0, 1.0 + rho, 0.5, 0.5 );

    return max( 0.0,
           min( 1.0,
           Phi2help( x, y, rho ) + Phi2help( y, x, rho ) ) );
}
\end{verbatim}
\end{small}
\hrule
\caption{C++ source code for evaluation of $\Phi_2(x,y;\varrho)$}
\label{source_Phi}
\end{figure}


\section{Discussion}

Evaluation of $\Phi_2(x,y;\varrho)$ as in Section \ref{sec_implement} will require (at most) four calls to an implementation of the cumulative standard normal distribution $\Phi$ ({\tt Phi()} in the code). The actual choice may well determine both accuracy and running time of the algorithm. For testing purposes I have been using a hybrid method, calling the algorithm from \cite[Fig. 2]{West} for absolute value larger than $0.5$, and {\tt Phi()} from \cite{Marsaglia} else. Besides {\tt Phi()}, {\tt exp()} will be called two times, {\tt arcsin()} or {\tt arccos()} two times, and {\tt sqrt()} six times. Everything else is elementary arithmetic.

Due to the reduction algorithm, the final result will be a sum. Therefore, very high accuracy in terms of relative error can not be expected. Consequently, evaluation of the diagonal aims at absolute error as well.

The {\tt Phi2diag()} function is behaving as it may be expected from an approximation by a Taylor series around zero: (absolute) error increases with decreasing $x$. For $x<-7$ (or $\varrho\longrightarrow 0$ or $\varrho\longrightarrow 1$) the error bounds from (\ref{eq_bounds}) are taking over, and absolute error decreases again. The maximum absolute error is obtained for $x\approx -7$, $\varrho\approx 0.8$ (maximum error of the upper bound is obtained for $\varrho=\sqrt{1-4/\pi^2}\approx 0.7712$, cf. \cite[Th. 5.2]{Meyer}).

In general, assuming that all numerical fallacies in the reduction algorithm have been taken care of, the diagonal is expected to provide a worst case because the errors of the two calls to {\tt Phi2diag()} will not cancel. With respect to the reduction algorithm, the case $\varrho x \approx y$, $\varrho y \approx x$, implying $|\varrho|\approx 1$, is most critical.

In order to give an impression of the algorithm's behaviour, we will discuss the results of a simulation study. For each $n\in\{0,\ldots,200\}$, $m\in\{1,\ldots,10^6\}$, the value of $\Phi_2(x_m^n,y_m^n;\varrho_m^n)$ has been computed via the {\tt Phi2()} function from Figure \ref{source_Phi} where $x_m^n$ has been drawn from a uniform distribution on $[x^n-0.05,x^n+0.05]$ with $x^n := n/10-10$, $y_m^n$ has been drawn from a uniform distribution on $[-10,10]$, and $\varrho_m^n:=2\Phi(r_m^n)-0.5$ where $r_m^n$ has been drawn from a uniform distribution on $[-10,10]$ as well.

The C++ implementation from \cite{West} has been serving as a competitor. Both functions have been evaluated against a quad-double precision version of {\tt Phi2()}, implemented using the QD library \cite{QD} and quad-double precision constants.

\begin{figure}
\begin{center}
\includegraphics[width=11.3cm]{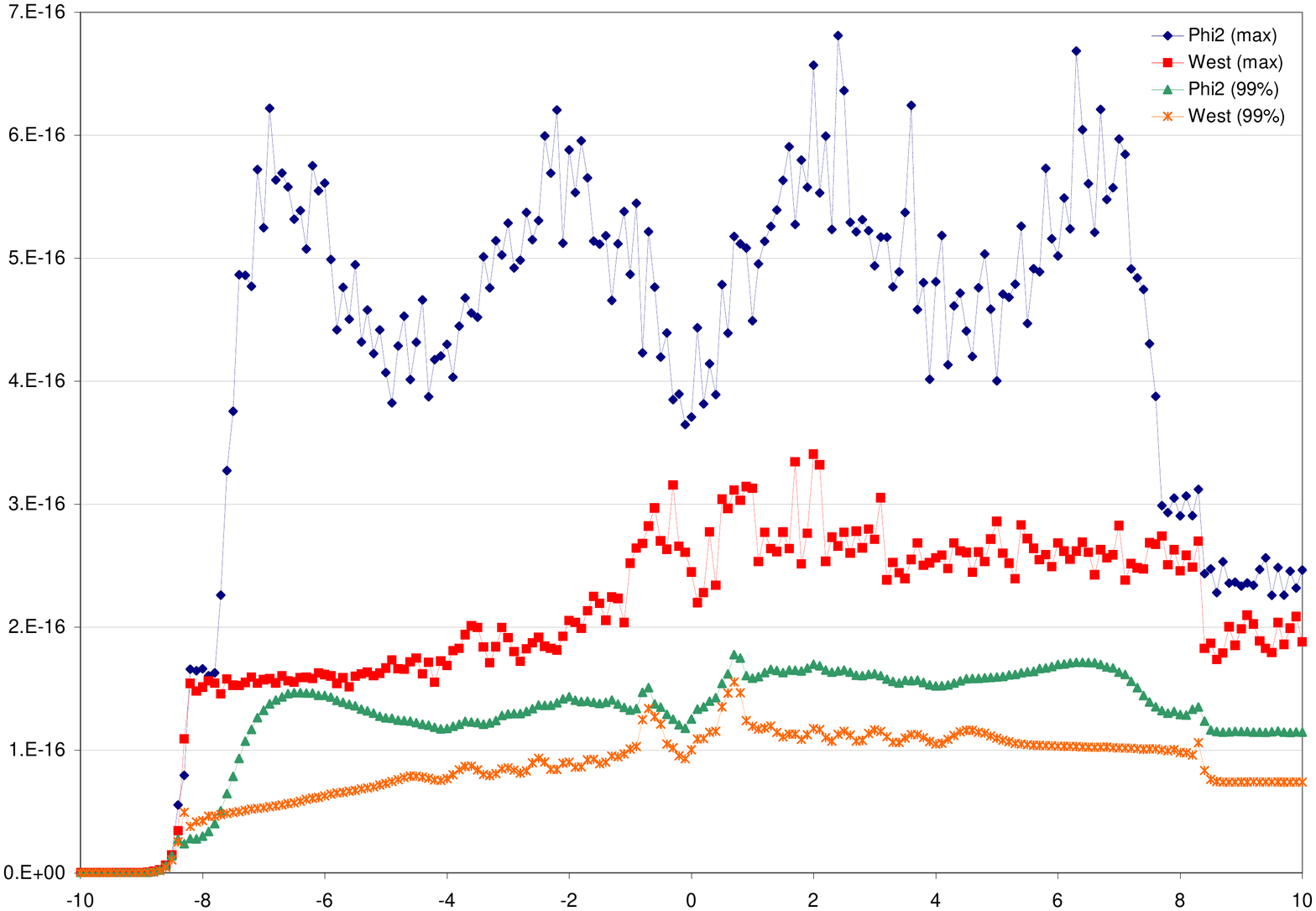}
\caption{Absolute error of implementations of $\Phi_2$ in a simulation study}
\label{diag_error}
\end{center}
\end{figure}

The diagram in Figure \ref{diag_error} is displaying, for $n\in\{0,\ldots,200\}$, the 99\% quantile and the maximum of the absolute difference between the double precision algorithms ({\tt Phi2} and {\tt West}) and the quad-double precision algorithm.

Apart from a shift due to subtractions for positive $x$, errors of {\tt Phi2} are rather symmetric around zero. The peaks at $|x^n|\approx 7$ are due to the Taylor expansion around zero; the peaks at $|x^n|\approx 2$ are due to Taylor expansion after transformation of the argument. The characteristics of the $99\%$ quantile, in particular the little peaks at $|x^n|\approx 0.7$, are already visible in the error of the $\Phi$ function used. The maximum error of {\tt West} almost always stays below the one of {\tt Phi2}. Note that the maximum error of {\tt West} is determined by the case $\varrho\longrightarrow -1$ and might be reduced by careful consideration of that case.

In the simulation study, {\tt Phi2} was a little slower than {\tt West}: it took approximately five minutes and four minutes to perform the $201\cdot 10^6$ evaluations on a fairly standard office PC (and it took two days to perform the corresponding quad-double precision evaluations). The number of recursion steps used by {\tt Phi2diag} is increasing with $|x|$. Because of the mathematical transparency of the algorithm it should be easy to find an appropriate trade-off between speed and accuracy by replacing the condition terminating the recursion.



\end{document}